# ON ISOPHOTE CURVE AND ITS CHARACTERIZATIONS

FATIH DOĞAN AND YUSUF YAYLI

ABSTRACT. Isophote comprises a locus of the surface points whose normal vectors make a constant angle with a fixed vector. Main objective of this paper is to find the axis of an isophote curve via its Darboux frame and afterwards to give some characterizations about isophote and its axis. Particularly, for isophotes lying on a canal surface will be obtained other characterizations again.

## 1. INTRODUCTION

Isophote is one of the characteristic curves on a surface such as parameter, geodesic and asymptotic curves or line of curvature.

Isophote on a surface can be regarded as a nice consequence of Lambert's cosine law in optics branch of physics. Lambert's law states that the intensity of illumination on a diffuse surface is proportional to the cosine of the angle generated between the surface normal vector $N$ and the light vector $d$. According to this law the intensity is irrespective of the actual viewpoint, hence the illumination is the same when viewed from any direction [15]. In other words, isophotes of a surface are curves with the property that their points have the same light intensity from a given source (a curve of constant illumination intensity). When the source light is at infinity, we may consider that the light flow consists in parallel lines. Hence, we can give a geometric description of isophotes on surfaces, namely they are curves such that the surface normal vectors in points of the curve make a constant angle with a fixed direction (which represents the light direction). These curves are succesfully used in computer graphics but also it is interesting to study for geometry. Then, to find an isophote on a surface we use the formula

$$\frac{\langle N(u,v), d \rangle}{\|N(u,v)\|} = \cos\theta, \ 0 \le \theta \le \frac{\pi}{2}.$$

In the special case, isophote is called as a silhouette curve if

$$\frac{\langle N(u,v), d \rangle}{\|N(u,v)\|} = \cos\frac{\pi}{2} = 0,$$

where $d$ is the vector of the line of sight from infinity.

Koenderink and van Doorn [7] studied the field of constant image brightness contours (isophotes). They showed that the spherical image (the Gauss map) of an isophote is a latitude circle on the unit sphere $S^2$ and the problem was reduced to that of obtaining the inverse Gauss map of these circles. By means of this they defined two kind singularities of the Gauss map: folds (curves) and simple cusps (apex, antapex points) and there are structural properties of the field of isophotes that bear an invariant relation to geometric features of the object.

Poeschl [12] used isophotes in car body construction via detecting irregularities along these curves on a free form surface. These irregularities (discontinuity of a surface or of the Gaussian curvature) emerge by differentiating of the equation $\langle N(u,v), l \rangle = \cos \theta = c$ (constant)

$$\langle N_u, l \rangle du + \langle N_v, l \rangle dv = 0$$

$$\frac{dv}{du} = -\frac{\langle N_u, l \rangle}{\langle N_v, l \rangle}, \ \langle N_v, l \rangle \neq 0,$$

where $l$ is the light vector.

Sara [14] researched local shading of a surface through isophotes properties. By using fundamental theory of surfaces, he focused on accurate estimation of surface normal tilt and on qualitatively correct Gaussian curvature recovery.

Kim and Lee [6] parameterized isophotes for surface of rotation and canal surface. They utilized that both these surfaces decompose into a set of circles where the surface normal vectors at points on each circle construct a cone. Again the vectors that make a constant angle with the fixed vector $d$ construct another cone and thus tangential intersection of these cones give the parametric range of the connected component isophote. In the same way, the same authors [5] parameterized the perspective silhouette of a canal surface by solving the problem that characteristic circles meet each other tangentially.

Dillen *et al.* [1] studied the constant angle surfaces in the product space $\mathbb{S}^2 \times \mathbb{R}$ for which the unit normal makes a constant angle with the $\mathbb{R}$-direction. Then Dillen and Munteanu [2] investigated the same problem in $\mathbb{H}^2 \times \mathbb{R}$ where $\mathbb{H}^2$ is the hyperbolic plane. Again, Nistor [10] researched normal, binormal and tangent developable surfaces of the space curve from standpoint of constant angle surface. Recently, Munteanu and Nistor [9] gave a important characterization about constant angle surfaces and studied the constant angle surfaces taking with a fixed vector direction being the tangent direction to $\mathbb{R}$ in Euclidean 3-space. Thus, it can be said that all curves on a constant angle surface are isophote curves.

Izumiya and Takeuchi [4] defined a slant helix as a space curve that the principal normal lines make a constant angle with a fixed direction. They displayed that a certain slant helix is also a geodesic on the tangent developable surface of a general helix. As an amazing consequence in our paper, we see that the curve which is both a geodesic and a slant helix on a surface is an isophote.

A canal surface is the envelope of a family of one parameter spheres and is useful to represent various objects e.g. pipe, hose, rope or intestine of a body. Canal surface is an important instrument in surface modelling for CAD/CAM such as tubular surfaces, torus and Dupin cyclides.

In this paper, we give some basic facts and concepts concerning curve and surface theory in section 2. In section 3, we concentrate on finding the axis of an isophote and also to characterize it in different ways. Here, we see that there is a close relation between isophotes and some special curves on a surface. Finally, in section 4 for an isophote lying on a canal surface, we obtain interesting results as to moving sphere which generates such a canal surface and then we find isophotes as some $v$-parameter curves on the tube.

## 2. Preliminaries

In the first place, we give some basic notions about curves and surfaces. The differential geometry of curves starts with a smooth map of $s$, let us call it $\alpha : I \subset \mathbb{R} \longrightarrow \mathbb{E}^3$ that parameterized a spatial curve denoted again with $\alpha$. We say that the curve is parameterized by arc-lenght if $\left\| \alpha'(s) \right\| = 1$ (unit speed), where $\alpha'(s)$ is the first derivative of $\alpha$ with respect to $s$. Let $\alpha : I \subset \mathbb{R} \longrightarrow \mathbb{E}^3$ be a regular curve with an arc-lenght parameter $s$ and $\kappa(s) = \left\| \alpha''(s) \right\| > 0$ where $\kappa$ is the curvature of $\alpha$ and $\alpha''$ is the second derivative of $\alpha$ with respect to $s$. Since the curvature $\kappa$ is nonzero, the Frenet frame $\{T, n, b\}$ is well-defined along the curve $\alpha$ and as follows.

$$\tag{2.1} T(s) = \alpha'(s)$$

$$n(s) = \frac{\alpha''(s)}{\|\alpha''(s)\|}$$

$$b(s) = T(s) \times n(s),$$

where $T$, $n$ and $b$ are the tangent, the principal normal and the binormal of $\alpha$, respectively. For a unit speed curve with $\kappa > 0$ the derivatives of the Frenet frame (Frenet-Serret formulas) are given by

$$\tag{2.2} T'(s) = \kappa(s) n(s)$$

$$n'(s) = -\kappa(s) t(s) + \tau(s) b(s)$$

$$b'(s) = -\tau(s) n(s),$$

where $\tau(s) = \dfrac{\left\langle \alpha'(s) \times \alpha''(s), \alpha'''(s) \right\rangle}{\kappa^2}$ is the torsion of $\alpha$; "$\times$" is the cross product on $\mathbb{R}^3$.

Let $M$ be a regular surface and $\alpha : I \subset \mathbb{R} \longrightarrow M$ be a unit speed curve on the surface. Then, the Darboux frame $\{T, B = N \times T, N\}$ is well-defined along the curve $\alpha$, where $T$ is the tangent of $\alpha$ and $N$ is the unit normal of $M$. Darboux equations of this frame are given by

$$\tag{2.3} T' = k_g B + k_n N$$

$$B' = -k_g T + \tau_g N$$

$$N' = -k_n T - \tau_g B,$$

where $k_n$, $k_g$ and $\tau_g$ are the normal curvature, the geodesic curvature and the geodesic torsion of $\alpha$, respectively. With the above notations, let $\phi$ denote the angle between the surface normal $N$ and the binormal $b$. Using equations in (2.3) we reach

$$\tag{2.4} \kappa^2 = k_g^2 + k_n^2$$

$$k_g = \kappa \cos \phi$$

$$k_n = \kappa \sin \phi$$

$$\tau_g = \tau - \phi'.$$

If we rotate the Darboux frame $\{T,\ B = N \times T,\ N\}$ by $\phi$ about $T$, we obtain the Frenet frame $\{T,\ n,\ b\}$.

$$\begin{bmatrix} T \\ n \\ b \end{bmatrix} = \begin{bmatrix} 1 & 0 & 0 \\ 0 & \cos\phi & \sin\phi \\ 0 & -\sin\phi & \cos\phi \end{bmatrix} \begin{bmatrix} T \\ B \\ N \end{bmatrix}$$

$$\begin{aligned} T &= T \\ n &= \cos\phi B + \sin\phi N \\ b &= -\sin\phi B + \cos\phi N. \end{aligned}$$

From the above equations, we obtain

(2.5)
$$N = \sin(\phi)n + \cos(\phi)b$$
$$B = \cos(\phi)n - \sin(\phi)b.$$

## 3. The Axis Of An Isophote Curve

In this section, we will get the fixed vector $d$ of an isophote curve via its Darboux frame. Let $M$ be a regular surface and let $\alpha : I \subset \mathbb{R} \longrightarrow M$ be a unit speed isophote curve. Then, from the definition of isophote curve

(3.1)
$$\langle N(u,v), d \rangle = \cos\theta = constant,$$

where $N(u,v)$ is the unit normal vector field of the surface $S(u,v)$ (a parameterization of $M$) and $d$ is the unit fixed vector on the axis of isophote curve.

Now, we begin to find the fixed vector $d$. Since $\alpha : I \subset \mathbb{R} \longrightarrow M$ is a unit speed isophote curve, the Darboux frame can be defined as $\{T,\ B = N \times T,\ N\}$ along the curve $\alpha$. If we differentiate Eq.(3.1) with respect to $s$ along the curve,

(3.2)
$$\left\langle N',d \right\rangle = 0.$$

From Eq.(2.3), it follows that

$$\langle -k_n T - \tau_g B, d \rangle = 0$$
$$-k_n \langle T,d \rangle - \tau_g \langle B,d \rangle = 0$$
$$\langle T,d \rangle = -\frac{\tau_g}{k_n} \langle B,d \rangle.$$

Because the Darboux frame $\{T,\ B,\ N\}$ is an orthonormal basis, if we say $\langle B,d \rangle = a$ in the last equation, then $d$ can be written as

$$d = -\frac{\tau_g}{k_n} aT + aB + \cos\theta N.$$

Since $\|d\| = 1$, we get

$$a = \mp \frac{k_n}{\sqrt{k_n^2 + \tau_g^2}} \sin\theta.$$

Thus, the vector $d$ is obtained as

(3.3)
$$d = \pm \frac{\tau_g}{\sqrt{k_n^2 + \tau_g^2}} \sin\theta T \mp \frac{k_n}{\sqrt{k_n^2 + \tau_g^2}} \sin\theta B + \cos\theta N$$

or from Eq.(2.4) and Eq.(2.5) in terms of the Frenet frame,

$$d = \pm \frac{\tau_g}{\sqrt{k_n^2 + \tau_g^2}} \sin\theta T + \left[\mp \frac{k_n k_g}{\kappa\sqrt{k_n^2 + \tau_g^2}} \sin\theta + \frac{k_n}{\kappa} \cos\theta\right] n$$

(3.4)
$$+ \left[\pm \frac{k_n^2}{\kappa\sqrt{k_n^2 + \tau_g^2}} \sin\theta + \frac{k_g}{\kappa} \cos\theta\right] b.$$

In fact, $d' = 0$ in other words $d$ is a constant vector. If we differentiate $N'$ and Eq.(3.2) with respect to $s$, we get

$$N'' = (-k_n' + k_g \tau_g)T - (k_n k_g + \tau_g')B - (k_n^2 + \tau_g^2)N$$

$$\left\langle N'', d \right\rangle = \frac{\mp(k_n' \tau_g - k_n \tau_g') \pm k_g(k_n^2 + \tau_g^2)}{\sqrt{k_n^2 + \tau_g^2}} \sin\theta - (k_n^2 + \tau_g^2)\cos\theta = 0,$$

where $d$ is form of Eq.(3.3). As a result, we have

$$\cot\theta = \pm \left[\frac{k_n^2}{(k_n^2 + \tau_g^2)^{\frac{3}{2}}} \left(\frac{\tau_g}{k_n}\right)' + \frac{k_g}{(k_n^2 + \tau_g^2)^{\frac{1}{2}}}\right]$$

(3.5)
$$\tan\theta = \frac{(k_n^2 + \tau_g^2)^{\frac{3}{2}}}{\pm k_g(k_n^2 + \tau_g^2) \pm (k_n \tau_g' - k_n' \tau_g)}.$$

By Eq.(2.3) and Eq.(3.3), the derivative of $d$ with respect to $s$ is that

$$d' = \pm \sin\theta \left[(\frac{\tau_g}{\sqrt{k_n^2 + \tau_g^2}})'T + \frac{\tau_g}{\sqrt{k_n^2 + \tau_g^2}}(k_g B + k_n N)\right]$$

$$\mp \sin\theta \left[(\frac{k_n}{\sqrt{k_n^2 + \tau_g^2}})'B + \frac{k_n}{\sqrt{k_n^2 + \tau_g^2}}(-k_g T + \tau_g N)\right] + \cos\theta(-k_n T - \tau_g B).$$

If we arrange this equality, we obtain

$$d' = \left(\pm \sin\theta \left[(\frac{\tau_g}{\sqrt{k_n^2 + \tau_g^2}})' + \frac{k_g k_n}{\sqrt{k_n^2 + \tau_g^2}}\right] - k_n \cos\theta\right) T$$

(3.6)
$$+ \left(\pm \sin\theta \left[-(\frac{k_n}{\sqrt{k_n^2 + \tau_g^2}})' + \frac{k_g \tau_g}{\sqrt{k_n^2 + \tau_g^2}}\right] - \tau_g \cos\theta\right) B.$$

Furthermore, from Eq.(3.5) we have

$$\cos\theta = \pm \sin\theta \frac{k_g(k_n^2 + \tau_g^2) + k_n \tau_g' - k_n' \tau_g}{(k_n^2 + \tau_g^2)^{\frac{3}{2}}}.$$

If the last equality is replaced in Eq.(3.6), we get

$$d' = \pm \sin\theta \left( \begin{array}{c} \dfrac{\tau'_g(k_n^2 + \tau_g^2) - \tau_g(k_n k'_n + \tau_g \tau'_g) + k_g k_n(k_n^2 + \tau_g^2)}{(k_n^2 + \tau_g^2)^{\frac{3}{2}}} \\ + \dfrac{k_n k'_n \tau_g - k_n^2 \tau'_g - k_g k_n(k_n^2 + \tau_g^2)}{(k_n^2 + \tau_g^2)^{\frac{3}{2}}} \end{array} \right) T$$

$$\pm \sin\theta \left( \begin{array}{c} \dfrac{-\kappa'_n(k_n^2 + \tau_g^2) + k_n(k_n k'_n + \tau_g \tau'_g) + k_g \tau_g(k_n^2 + \tau_g^2)}{(k_n^2 + \tau_g^2)^{\frac{3}{2}}} \\ + \dfrac{k'_n \tau_g^2 - k_n \tau_g \tau'_g - k_g \tau_g(k_n^2 + \tau_g^2)}{(k_n^2 + \tau_g^2)^{\frac{3}{2}}} \end{array} \right) B$$

As can be directly seen above, the coefficients of $T$ and $B$ becomes zero. So, $d' = 0$ namely $d$ is a constant vector. Then, the axis of isophote is the line in fixed direction $d$. From this time, for the axis of isophote will be used $d$ again.

**Theorem 1.** *A unit speed curve $\alpha$ on a surface is an isophote if and only if*

$$\cot\theta = \mu(s) = \pm \left( \dfrac{k_n^2}{(k_n^2 + \tau_g^2)^{\frac{3}{2}}} \left(\dfrac{\tau_g}{k_n}\right)' + \dfrac{k_g}{(k_n^2 + \tau_g^2)^{\frac{1}{2}}} \right)(s)$$

*is a constant function.*

*Proof.* As $\alpha$ is an isophote, the Gauss map along the curve $\alpha$ is a circle on the unit sphere $S^2$. Hence, if we compute the Gauss map $N_{|\alpha} : I \longrightarrow S^2$ along the curve $\alpha$, the geodesic curvature of $N_{|\alpha}$ becomes $\mu(s)$ as shown below.

$$\begin{aligned} N'_{|\alpha} &= -k_n T - \tau_g B \\ N''_{|\alpha} &= (-k'_n + k_g \tau_g)T - (k_n k_g + \tau'_g)B - (k_n^2 + \tau_g^2)N \\ N'_{|\alpha} \times N''_{|\alpha} &= \tau_g(k_n^2 + \tau_g^2)T + \left(k_g(k_n^2 + \tau_g^2) + k_n^2\left(\dfrac{\tau_g}{k_n}\right)'\right)N - k_n(k_n^2 + \tau_g^2)B \end{aligned}$$

Therefore, we get

$$\begin{aligned} \kappa &= \dfrac{\left\| N'_{|\alpha} \times N''_{|\alpha} \right\|}{\left\| N'_{|\alpha} \right\|^3} \\ &= \dfrac{\sqrt{\tau_g^2(k_n^2 + \tau_g^2)^2 + \left(k_g(k_n^2 + \tau_g^2) + k_n^2\left(\dfrac{\tau_g}{k_n}\right)'\right)^2 + k_n^2(k_n^2 + \tau_g^2)^2}}{\sqrt{(k_n^2 + \tau_g^2)^3}} \\ &= \sqrt{1 + \dfrac{\left(k_g(k_n^2 + \tau_g^2) + k_n^2\left(\dfrac{\tau_g}{k_n}\right)'\right)^2}{(k_n^2 + \tau_g^2)^3}}. \end{aligned}$$

Let $\bar{k}_g$ and $\bar{k}_n$ be the geodesic curvature and the normal curvature of the Gauss map $N_{|\alpha}$, respectively. Since the normal curvature $\bar{k}_n = 1$ on $S^2$, if we substitute $\bar{k}_n$ and $\kappa$ in the following equation, we obtain the geodesic curvature $\bar{k}_g$ as follows.

$$\kappa^2 = (\bar{k}_g)^2 + (\bar{k}_n)^2$$

$$\bar{k}_g(s) = \mu(s) = \cot\theta = \pm\left(\frac{k_n^2}{(k_n^2+\tau_g^2)^{\frac{3}{2}}}(\frac{\tau_g}{k_n})' + \frac{k_g}{(k_n^2+\tau_g^2)^{\frac{1}{2}}}\right)(s).$$

Then, the spherical image of isophotes are circles if and only if $\mu(s)$ is a constant function. □

**Lemma 1.** *Let $\alpha$ be a unit speed space curve with $\kappa(s) \neq 0$. Then, $\alpha$ is a slant helix if and only if $\sigma(s) = \left(\frac{\kappa^2}{(\kappa^2+\tau^2)^{\frac{3}{2}}}(\frac{\tau}{\kappa})'\right)(s)$ is a constant function.*[4]

**Lemma 2.** *Let $\alpha$ be a unit speed space curve with $\kappa(s) \neq 0$. Then, $\alpha$ is a general helix if and only if $\frac{\tau}{\kappa}(s)$ is a constant function.*

**Theorem 2.** *Let $\alpha$ be a unit speed isophote curve on the surface $M$. In that case, we have the following.*
*(1) $\alpha$ is a geodesic on $M$ if and only if $\alpha$ is a slant helix with the fixed vector*

$$d = \pm\frac{\tau}{\sqrt{\kappa^2+\tau^2}}\sin\theta T \pm \cos\theta n \pm \frac{\kappa}{\sqrt{\kappa^2+\tau^2}}\sin\theta b.$$

*(2) $\alpha$ is a asymptotic curve on $M$ if and only if $\alpha$ is a general helix with the fixed vector $d = \pm\sin\theta T \pm \cos\theta b$.*
*(3) If $\alpha$ is a line of curvature, then $\alpha$ is a plane curve and the angle $\theta = \mp\phi$ or $\theta = \mp(\pi-\phi)$.*

*Proof.* (1) Since $\alpha$ is a geodesic (i.e. the surface normal $N$ concurs with the principal normal $n$ along the curve $\alpha$), we have $k_g = 0$ and therefore from Eq.(2.4) $k_n = \pm\kappa$ and $\tau_g = \tau$. By substituting $k_g$ and $k_n$ in the expression of $\mu(s)$, we follow that

$$\mu(s) = \pm\left(\frac{\kappa^2}{(\kappa^2+\tau^2)^{\frac{3}{2}}}(\frac{\tau}{\kappa})'\right)(s)$$

is a constant function. Then, by Lemma (1) $\alpha$ is a slant helix. Using Eq.(3.4), the fixed vector of the slant helix is obtained as

$$d = \pm\frac{\tau}{\sqrt{\kappa^2+\tau^2}}\sin\theta T \pm \cos\theta n \pm \frac{\kappa}{\sqrt{\kappa^2+\tau^2}}\sin\theta b.$$

On the contrary, let $\alpha$ be a slant helix with the fixed vector

$$d = \pm\frac{\tau}{\sqrt{\kappa^2+\tau^2}}\sin\theta T \pm \cos\theta n \pm \frac{\kappa}{\sqrt{\kappa^2+\tau^2}}\sin\theta b.$$

In this case, from Eq.(3.4) again the geodesic curvature $k_g$ must be zero that is to say $\alpha$ is a geodesic on $M$.

(2) Since $\alpha$ is a asymptotic curve on $M$, we have $k_n = 0$ and consequently from Eq.(2.4) $k_g = \pm\kappa$ and $\tau_g = \tau$. If we replace $k_n$, $k_g$ and $\tau_g$ in Eq.(3.5), we obtain that

$$\tan\theta = \pm\frac{\tau_g}{k_g} = \pm\frac{\tau}{\kappa}(s)$$

is a constant function. Then, from Lemma (2) and Eq.(3.4) it follows that $\alpha$ is a general helix with the fixed vector

$$d = \pm \sin\theta T + \cos\theta b.$$

For the converse, let $\alpha$ be a general helix with the fixed vector $d = \pm \sin\theta T \pm \cos\theta b$. By applying Eq.(3.4) again, we get $k_n = 0$ in other words $\alpha$ is a asymptotic curve on $M$.

(3) Since $\alpha$ is a line of curvature on $M$, we have $\tau_g = 0$. Accordingly, from Eq.(3.5) we conclude that

$$\tan\theta = \pm\frac{k_n}{k_g} = \pm\frac{\kappa\sin\phi}{\kappa\cos\phi} = \pm\tan\phi.$$

In this situation, we trace that $\phi = \pm\theta$ or $\phi = \pi \pm \theta$. Because $\phi$ is a constant, by $\tau_g = \tau - \phi' = 0$, we reach $\tau = 0$. Then, $\alpha$ is a plane curve. □

We now illustrate Theorem (2) case (1) in the following example.

**Example 1** ([4]). *We consider a space curve defined by*

$$\gamma(\theta) = \left(-\frac{(a^2-b^2)}{2a}\left(\frac{\cos((a+b)\theta)}{(a+b)^2}+\frac{\cos((a-b)\theta)}{(a-b)^2}\right),\right.$$
$$\left.-\frac{(a^2-b^2)}{2a}\left(\frac{\sin((a+b)\theta)}{(a+b)^2}+\frac{\sin((a-b)\theta)}{(a-b)^2}\right),-\frac{\sqrt{a^2-b^2}}{ab}\cos(b\theta)\right).$$

*We can calculate that*

$$\kappa(\theta) = \sqrt{a^2-b^2}\cos(b\theta), \quad \tau(\theta) = \sqrt{a^2-b^2}\sin(b\theta)$$

$$\sigma(\theta) = \frac{b}{\sqrt{a^2-b^2}}, \quad \frac{\tau}{\kappa}(\theta) = \tan(b\theta)$$

$$\left(\frac{\tau}{\kappa}(\theta)\right)' = \frac{b}{\cos^2(b\theta)} \neq 0, \text{ and } \left(\frac{\tau}{\kappa}(\theta)\right)'' = \frac{2b^2\tan(b\theta)}{\cos^2(b\theta)} \neq 0.$$

*Therefore, $\gamma(\theta)$ is a slant helix and it is not a cylindrical helix. By Theorem 4.2 [4], it is a geodesic of the tangent developable surface of a cylindrical helix. In fact, the corresponding tangent developable surface is the rectifying developable surface of $\gamma(\theta)$ by Proposition 4.1 [4]. We now draw the picture of $\gamma(\theta)$ ($a = 2, b = 1$) in Fig. 1a. We also draw the rectifying developable surface of $\gamma(\theta)$ in Fig. 1b.*

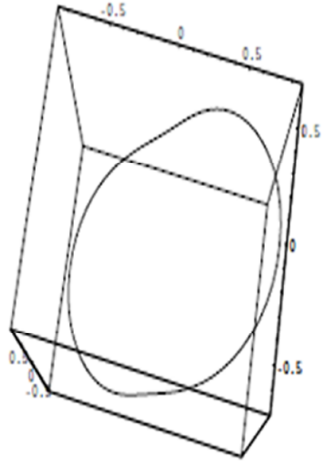

**Figure 1a**. The slant helix $\gamma(\theta)$ which is also a geodesic (Isophote) of the tangent developable surface of a cylindrical helix

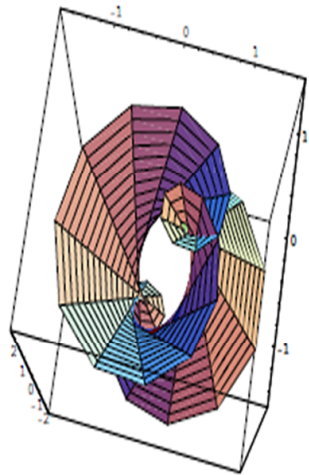

**Figure 1b.** The rectifying developable surface of $\gamma(\theta)$ corresponds to tangent developable surface of a cylindrical helix

In that case, we can say that there is a curve which is both a slant helix and a geodesic on a surface. This example is a consequence of Theorem (2) case (1) that is to say such a curve on a surface is an isophote. As a simpler example, the helix on a circular cylinder is a geodesic. Moreover, each helix is also a slant helix. Thus, the helix on a circular cylinder is both a slant helix and a geodesic to put a different way it is an isophote.

**Theorem 3.** *Let $\alpha$ be a unit speed isophote curve on the surface $M$. Then, we have the following.*
*(1) The axis $d$ is perpendicular to the tangent line of $\alpha$ if and only if $\alpha$ is a line of curvature on $M$.*

(2) The axis d is perpendicular to the principal normal line of $\alpha$ if and only if $\alpha$ is a asymptotic curve on $M$ or $\dfrac{\tau_g}{k_n}$ is a constant function.

*Proof.* (1) Let $\alpha$ be a unit speed isophote curve. Then by Eq.(3.4), it follows that

$$\langle T, d \rangle = \pm \dfrac{\tau_g}{\sqrt{k_n^2 + \tau_g^2}} \sin\theta.$$

Therefore, the axis $d$ is perpendicular to the tangent line of $\alpha$ if and only if $\alpha$ is a line of curvature on $M$.

(2) Let the axis $d$ be perpendicular to the principal normal line of $\alpha$. Then by Eq.(3.4), we have

$$\begin{aligned}\langle n, d \rangle &= \mp \dfrac{k_n k_g}{\kappa \sqrt{k_n^2 + \tau_g^2}} \sin\theta + \dfrac{k_n}{\kappa} \cos\theta \\ &= \dfrac{k_n}{\kappa} \left( \mp \dfrac{k_g}{\sqrt{k_n^2 + \tau_g^2}} \sin\theta + \cos\theta \right) \\ &= 0.\end{aligned}$$

By solving this equation, we find $k_n = 0$ or $\cot\theta = \pm \dfrac{k_g}{\sqrt{k_n^2 + \tau_g^2}}$. From Eq.(3.5), we gather that

$$\cot\theta = \dfrac{k_g}{\sqrt{k_n^2 + \tau_g^2}} = \dfrac{k_n^2}{(k_n^2 + \tau_g^2)^{\frac{3}{2}}} \left(\dfrac{\tau_g}{k_n}\right)' + \dfrac{k_g}{(k_n^2 + \tau_g^2)^{\frac{1}{2}}}$$

$$\dfrac{k_n^2}{(k_n^2 + \tau_g^2)^{\frac{3}{2}}} \left(\dfrac{\tau_g}{k_n}\right)' = 0.$$

In the last equation, for $k_n \neq 0$, $\left(\dfrac{\tau_g}{k_n}\right)'(s) = 0$ and hence $\left(\dfrac{\tau_g}{k_n}\right)(s)$ is a constant function. In the same way, the proof of sufficiency is clear. □

**Corollary 1.** *Let $\alpha$ be a silhouette curve on $M$. Then, $\alpha$ is a line of curvature if and only if it is a plane geodesic curve.*

*Proof.* Let $\alpha$ be a line of curvature. Because $\alpha$ is both a silhouette curve and a line of curvature, we possess $\tau_g = 0$ and $\phi = \theta = \pm\dfrac{\pi}{2}$. Therefore, $k_g = 0$ namely $\alpha$ is a geodesic. Now that $\phi$ is a constant and $\tau_g = 0$, from Eq.(2.4) $\tau$ must be zero. Eventually, $\alpha$ is a plane geodesic curve.

Conversely, let $\alpha$ be a plane geodesic curve. In this case, $k_g$ and $\tau$ need to be zero. From this, $\tau_g = 0$ in other words $\alpha$ is a line of curvature. This completes the proof. Furthermore, in the above corollary by Eq.(2.3) it follows that $B' = 0$, i.e., by Eq.(3.3) the axis of silhouette $d$ becomes $B$. □

**Corollary 2.** *Let $\alpha$ be a silhouette with arc-lenght parameter on $M$. In this case,*
*(1) The axis $d$ lies in the plane spanned by $T$ and $B$.*
*(2) $\alpha$ is a geodesic if and only if the axis $d$ lies on the rectifying plane of the silhouette curve $\alpha$.*

*Proof.* (1) Since $\alpha$ is a silhouette curve, the surface normal vectors are orthogonal to the axis $d$ that is $\theta = \dfrac{\pi}{2}$. Then, by Eq.(3.3) it follows that

$$d = \pm \frac{\tau_g}{\sqrt{k_n^2 + \tau_g^2}} T \mp \frac{k_n}{\sqrt{k_n^2 + \tau_g^2}} B.$$

We see that $d$ lies in the plane spanned by $T$ and $B$.

(2) Since $\alpha$ is a geodesic, by Theorem (2) $\alpha$ is a slant helix with the fixed vector

$$d = \pm \frac{\tau}{\sqrt{\kappa^2 + \tau^2}} \sin\theta T + \cos\theta n \mp \frac{\kappa}{\sqrt{\kappa^2 + \tau^2}} \sin\theta b.$$

Besides, as $\alpha$ is a silhouette curve, we have $\theta = \dfrac{\pi}{2}$ and

$$d = \pm \frac{\tau}{\sqrt{\kappa^2 + \tau^2}} T \mp \frac{\kappa}{\sqrt{\kappa^2 + \tau^2}} b$$

in other words the axis $d$ lies on the rectifying plane of the silhouette curve $\alpha$. By contrast, suppose that the axis $d$ lies on the rectifying plane of the silhouette curve $\alpha$. For this reason, applying Eq.(3.4) we see that $k_g = 0$ namely $\alpha$ is a geodesic on $M$. □

## 4. Some Characterizations for Isophotes on a Canal Surface

In this section, we introduce canal surfaces and tubes and then give some characterizations for isophotes on them. To begin with, we define a canal surface. Canal surface is defined as the envelope of a family of one parameter spheres. Alternatively, a canal surface is the envelope of a moving sphere with varying radius, defined by the trajectory $C(t)$ (spine curve) of its centers and a radius function $r(t)$. If the radius function $r(t) = r$ is a constant, then the canal surface is called a tube or pipe surface.

Now, we give parametric representation of a canal surface. Since canal surface is the envelope of one parameter spheres with center $C(t)$ and radius $r(t)$, a surface point $p \in \mathbb{E}^3$ satisfies the following equations.

$$\|p - C(t)\| = r(t)$$

(4.1) $$(p - C(t)) \cdot C'(t) + r(t)r'(t) = 0.$$

If the spine curve $C(t)$ has arc-lenght parametrization, then the canal surface is parametrized as

(4.2) $$K(s,v) = C(s) - r(s)r'(s)T(s) \mp r(s)\sqrt{1 - r'(s)^2}\left(\cos v\, n(s) + \sin v\, b(s)\right),$$

where $0 \leq v < 2\pi$; $n$ and $b$ are the principal normal and the binormal to $C(t)$, respectively.

Since the moving sphere with center $C(s)$ tangent to canal surface, unit normal vector field of the canal surface can be computed as follows.

$$N(s,v) = K(s,v) - C(s)$$

$$N(s,v) = -r(s)r'(s)T(s) \mp r(s)\sqrt{1 - r'(s)^2}\left(\cos v\, n(s) + \sin v\, b(s)\right).$$

As $\|K(s,v) - C(s)\| = r(s)$, the norm of the vector $N(s,v)$ is $r(s)$. Hence, if we normalize $N(s,v)$, we get

$$\text{(4.3)} \qquad \frac{N(s,v)}{\|N(s,v)\|} = -r'(s)T(s) \mp \sqrt{1 - r'(s)^2}\left(\cos v\, n(s) + \sin v\, b(s)\right).$$

In the previous section, we mentioned the axis of isophotes and obtained some characterizations of them. This time, we shall give other characterizations as regards isophotes lying on a canal surface, particularly.

**Theorem 4.** *If $\alpha$ is a unit speed isophote with the axis $d$ on the canal surface $K(s,v)$, then*

$$-r'\langle T,d\rangle \mp \sqrt{1-r'^2}\cos(v+\phi)\langle B,d\rangle + \left(\sqrt{1-r'^2}\sin(v+\phi) - 1\right)\langle N,d\rangle = 0.$$

*where $\phi$ is the angle between the surface normal $N$ and the binormal $b$.*

*Proof.* Inner product of the axis $d$ and the normal vector field $N$ is that

$$\langle N,d\rangle = -r'\langle T,d\rangle \mp \sqrt{1-r'^2}\left[\cos v\,\langle n,d\rangle + \sin v\,\langle b,d\rangle\right].$$

If we substitute $\langle n,d\rangle$ and $\langle b,d\rangle$ above, we get

$$\langle N,d\rangle = -r'\langle T,d\rangle \mp \sqrt{1-r'^2}\left[\mp\frac{k_n k_g}{\kappa\sqrt{k_n^2 + \tau_g^2}}\sin\theta + \frac{k_n}{\kappa}\cos\theta\right]\cos v$$

$$\mp \sqrt{1-r'^2}\left[\pm\frac{k_n^2}{\kappa\sqrt{k_n^2+\tau_g^2}}\sin\theta + \frac{k_g}{\kappa}\cos\theta\right]\sin v.$$

By arranging this equation, we find

$$\langle N,d\rangle = -r'\langle T,d\rangle \mp \sqrt{1-r'^2}\left[\frac{k_n\,(k_g\cos v - k_n\sin v)}{\kappa\sqrt{k_n^2+\tau_g^2}}\right]\sin\theta$$

$$+ \sqrt{1-r'^2}\left[\frac{k_n\cos v + k_g\sin v}{\kappa}\right]\cos\theta.$$

If we substitute $k_g = \kappa\cos\phi$, $k_n = \kappa\sin\phi$, $\langle B,d\rangle = \dfrac{k_n}{\sqrt{k_n^2+\tau_g^2}}$, $\langle N,d\rangle = \cos\theta$ and use addition formulas for sine and cosine, we obtain

$$-r'\langle T,d\rangle \mp \sqrt{1-r'^2}\cos(v+\phi)\langle B,d\rangle + \left(\sqrt{1-r'^2}\sin(v+\phi) - 1\right)\langle N,d\rangle = 0.$$

□

**Corollary 3.** *Let $\alpha$ be a unit speed isophote curve on the canal surface $K(s,v)$. Then we have the following.*
*(1) If the axis $d$ is orthogonal to the tangent line of $\alpha$, then the canal surface is generated by a moving sphere with linear radius function $r(s) = \lambda s + c$ where,*

$$\lambda = \frac{\sqrt{(\tan\theta\cos(v+\theta) \mp \sin(v+\theta))^2 - 1}}{\tan\theta\cos(v+\theta) \mp \sin(v+\theta)} \quad \text{and}$$

$$|\sin(v+2\theta)| > \cos\theta, \quad \cos\theta < \sin v < -\cos\theta.$$

(2) If $\alpha$ is a silhouette curve and the spine curve $C(s)$ is a general helix, then the canal surface is generated by a moving sphere with the radius function
$$r(s) = \int \frac{\tan\beta}{\sqrt{\tan^2\beta + \cos^2(v+\phi)}} ds + c; \ c > 0.$$

*Proof.* (1) Assume that the axis $d$ is orthogonal to the tangent line of $\alpha$ namely $d$ lies in the plane spanned by $N$ and $B$. Because $\langle N, d \rangle = \cos\theta$, by the preceding theorem we get
$$\sqrt{1-r'^2}\cos(v+\phi)\langle B,d\rangle \mp \left(\sqrt{1-r'^2}\sin(v+\phi)-1\right)\langle N,d\rangle = 0$$
$$\sqrt{1-r'^2}\cos(v+\phi)\sin\theta \mp \left(\sqrt{1-r'^2}\sin(v+\phi)-1\right)\cos\theta = 0.$$

Since $\langle T, d \rangle = 0$, from Theorem (2) and Theorem (3), $\alpha$ is a line of curvature with $\phi = \theta$ (constant). So, if the final equation is arranged, it follows that
$$1 - r'^2 = \frac{1}{\cos^2(v+\theta)(\tan\theta \mp \tan(v+\theta))^2}.$$

If we solve this quadratic equation with unknown $r'$, we obtain
$$r(s) = \left(\frac{\sqrt{(\tan\theta\cos(v+\theta)\mp\sin(v+\theta))^2 - 1}}{\tan\theta\cos(v+\theta)\mp\sin(v+\theta)}\right)s + c; \ c > 0.$$

For $(\tan\theta\cos(v+\theta)\mp\sin(v+\theta))^2 - 1 > 0$, $r(s) > 0$. By solving this inequality, we can reach the condition $|\sin(v+2\theta)| > \cos\theta$ and $\cos\theta < \sin v < -\cos\theta$.

(2) Assume that $\alpha$ is a silhouette curve namely $\langle N, d\rangle = 0$. Additionally, since $C(s)$ is a general helix, we can write $\langle T, d\rangle = \cos\beta$ ($\beta$ constant). Then
$$-r'\langle T,d\rangle \mp \sqrt{1-r'^2}\cos(v+\phi)\langle B,d\rangle = 0$$
$$-r'\cos\beta \mp \sqrt{1-r'^2}\cos(v+\phi)\sin\beta = 0.$$

If the last equation is arranged, the solution of the quadratic equation with unknown $r'$ is obtained as follows.
$$\left(\tan^2\beta + \cos^2(v+\phi)\right)r'^2 - \tan^2\beta = 0$$
$$r(s) = \int \frac{\tan\beta}{\sqrt{\tan^2\beta + \cos^2(v+\phi)}} ds + c; \ c > 0.$$

Since $\beta$ is an acute angle, $\tan\beta > 0$. Then, $r(s) > 0$. $\square$

**Proposition 1.** *Let the spine curve $C(s)$ be a general helix. If an isophote on the canal surface and the general helix $C(s)$ have the same axis $d$, then the canal surface is generated by a moving sphere with linear radius function $r(s) = \omega s + c$ where $\omega = \dfrac{-1 + \sin^2 v \tan\theta}{1 + \sin^2 v \tan^2\theta}$ and $\tan\theta > 1$.*

*Proof.* In that $C(s)$ is a general helix with the axis $d$, $\langle T, d\rangle = \cos\theta$ (constant) and so $\langle n, d\rangle = 0$. Thus
$$\langle N, d\rangle = \langle T, d\rangle = \cos\theta = \lambda_1, \ \langle b, d\rangle = \sin\theta = \lambda_2.$$

If we put $\lambda_1$ and $\lambda_2$ in Eq.(4.3), we find
$$\left(\lambda_1^2 + \lambda_2^2 \sin^2 v\right) r'^2 + 2\lambda_1^2 r' + \lambda_1^2 - \lambda_2^2 \sin^2 v = 0.$$

If we solve this quadratic equation, we obtain

$$r(s) = \left(\frac{-1 + \sin^2 v \tan \theta}{1 + \sin^2 v \tan^2 \theta}\right) s + c; \ c > 0.$$

Since $\theta$ is an acute angle, $\tan \theta > 0$. In addition, for $-1 + \sin^2 v \tan \theta > 0$, $r(s) > 0$. Therefore, it must be $\tan \theta > 1$.

From now on, we will give some characterizations for isophote on a tube. If the radius function $r(s) = r$ is a constant, by Eq.(4.2) a tube is parametrized as follows.

$$K(s,v) = C(s) \mp r \left(\cos v \ n(s) + \sin v \ b(s)\right).$$

□

**Proposition 2.** *Let the spine curve $C(s)$ be a general helix with the axis $d$. Then $v_0 = \left(\dfrac{2k+1}{2}\right)\pi$ ($k \in Z$) parameter curves of tube are isophotes with the axis $d$.*

*Proof.* Because the normal vector field of tube $N(s,v) = K(s,v) - C(s)$, we obtain $N(s,v) = \mp r \left(\cos v \ n(s) + \sin v \ b(s)\right)$. For $v_0 = \left(\dfrac{2k+1}{2}\right)\pi$,

$$N(s,v_0) = \mp r \ b(s).$$

Since $C(s)$ is a general helix with the axis $d$, $\langle b, d \rangle$ is a constant. Therefore, along the curve $v_0 = \left(\dfrac{2k+1}{2}\right)\pi$

$$\langle N(s,v_0), d \rangle = \mp r \langle b, d \rangle = constant.$$

In this situation $v_0 = \left(\dfrac{2k+1}{2}\right)\pi$ parameter curves are isophotes with the axis $d$ on the tube.

**Proposition 3.** *Let the spine curve $C(s)$ be a slant helix with the axis $d$. Then $v_0 = k\pi$ ($k \in Z$) parameter curves of tube are isophotes with the axis $d$.*

□

*Proof.* Because the normal vector field of tube $N(s,v) = K(s,v) - C(s)$, we obtain $N(s,v) = \mp r \left(\cos v \ n(s) + \sin v \ b(s)\right)$. For $v_0 = k\pi$,

$$N(s,v_0) = \mp r \ n(s).$$

Since $C(s)$ is a slant helix with the axis $d$, $\langle n, d \rangle$ is a constant. Therefore, along the curve $v_0 = k\pi$

$$\langle N(s,v_0), d \rangle = \mp r \langle n, d \rangle = constant.$$

Then $v_0 = k\pi$ parameter curves are isophotes with the axis $d$. □

## 5. Conclusions

In this paper, we found the axis (fixed vector) of an isophote curve through its Darboux frame. Subsequently, we obtained some characterizations regarding these curves. By using the characterizations, it is investigated relation between special curves on a surface and an isophote. Finally, we reached some results for an isophote lying on canal surface and then obtained several isophotes as special parameter curves on the tube.

*Current Adress: Fatih DOĞAN, Yusuf YAYLI, Ankara University, Department of Mathematics, 06100 Tandoğan, Ankara, Turkey*

*E-mail address*: mathfdogan@hotmail.com, yayli@science.ankara.edu.tr